\documentclass{amsart}
\usepackage{amssymb}
\usepackage{amscd}

\usepackage[
hyperref
]{degt}

\let\leftbracket[
{\catcode`\[\active\catcode`\*\active
\gdef\deMaple{\catcode`\[\active\catcode`\*\active
 \def[##1]{_{##1}}\let*\relax\let\[\leftbracket}}
 
\def\2{^{\prime2}}
\def\={\noalign{\vspace\medskipamount}}
\def\bx{\bar x}
\def\by{\bar y}
\def\bz{\bar z}
\def\tr{\tilde r}
\def\tC{\tilde C}
\def\tK{\tilde K}
\def\tQ{\tilde Q}
\def\tR{\tilde R}
\def\tS{\tilde S}
\let\da\dashrightarrow

\title{The Bertini involution}

\author{Alex Degtyarev}

\address{%
Bilkent University\\
Department of Mathematics\\
06800 Ankara, Turkey}

\email{degt@fen.bilkent.edu.tr}

\keywords{Bertini involution, Geiser involution}

\subjclass[2000]{14E07}

\thanks{These notes are to be extended should there be any interesting
development. They will be available at
\url{http://www.fen.bilkent.edu.tr/~degt/papers/papers.htm}
and on the {\tt arXiv}}

\begin{document}

\begin{abstract}
We summarize and extend E.~Moody's results on the explicit equations related
to the Bertini involution.
\end{abstract}

\maketitle

These notes are the result of my attempt to understand E.~Moody's
paper~\cite{Moody}. I correct a few misprints in~\cite{Moody} and take the
computation a bit further.

I express my admiration to Ethel I.~Moody, who managed to perform this
tedious computation in the pre-\Maple\ era.
A \Maple\ implementation
of most equations is found at
\url{http://www.fen.bilkent.edu.tr/~degt/papers/Bertini.zip}.

This text is not intended as an `official' publication; it is distributed in
the hope that it may be useful. It can be cited by its {\tt arXiv} location.

Whenever possible, I try to keep the
original notation of~\cite{Moody}.

\section{The Bertini involution}\label{S.Bertini}

\subsection{The results of~\cite{Moody}}\label{s.Moody}
Consider the pencil of cubics
\[
\Gl w(x)+\mu w'(x)=0,
\label{eq.pencil}
\]
where
\[*
w(x)=x_3^2(a_1x_1+a_2x_2)+
x_3(b_1x_1^2+b_2x_1x_2+b_3x_2^2)+(c_1x_1^2x_2+c_2x_1x_2^2)
\]
and similar for $w'$, so that the coordinate vertices are amongst the
basepoints of the pencil.
The point $(0:0:1)$ will play a special r\^{o}le.

The curve of the pencil passing through a point~$y$
is given by
\[
W_3(x):=w(x)w'(y)-w'(x)w(y)=0.
\label{eq.W}
\]
Clearly,
\[*
W_3(x)=x_3^2(A_1x_1+A_2x_2)+
x_3(B_1x_1^2+B_2x_1x_2+B_3x_2^2)+(C_1x_1^2x_2+C_2x_1x_2^2),
\]
where $A_i(y):=a_iw'(y)-a'_iw(y)$ and similar for $B_i$, $C_i$.

The tangent to~\eqref{eq.W} at $(0:0:1)$ meets the curve again at
$r=(r_1:r_2:r_3)$, where
\[
\gathered
r_1=A_2r_1',\quad r_2=-A_1r_1',\quad r_3=A_1A_2r_3',\\
r_1':=B_1A_2^2-B_2A_1A_2+B_3A_1^2,\qquad
r_3':=A_2C_1-A_1C_2.
\endgathered
\label{eq.R}
\]
The locus of these points is
\[
\Gg_4(y):=y_1A_1+y_2A_2=0.
\label{eq.gamma}
\]
Apart from the basepoints, the locus~\eqref{eq.gamma} meets~\eqref{eq.W} at a
single point~$r$. The line $(ry)$ meets~\eqref{eq.W} at a third point~$z$, and
the \emph{Bertini involution} can be defined as the map $y\mapsto z$.
Let $\kappa:=a_1b'_1 - a'_1b_1$ and
\[*
\deMaple
\aligned
C_5(y)&:=A[2]*\[B[1] + \kappa*y[1]*y[3]^2]_{y[2]} +
  \[A[1] - \kappa*y[1]^2*y[3]]_{y[2]}*\[A[2]*y[3]+B[3]*y[2]]_{y[1]} +
  \kappa*B[3]*y[1]*y[3],\\
\phi_6(y)&:=A_1C_2+y[3]C_5(y),\\
\psi_6(y)&:=A_2C_1+y[3]C_5(y).
\endaligned
\]
(Following~\cite{Moody}, we use $[e]_u$ do indicate that $e$
has a common factor~$u$ and this factor has been removed.)
In these notations, the Bertini involution is
\[
\deMaple
 z[1] = \phi_6*\[A[2]^2*\phi_6 + B[3]*r'[1]]_{y[1]},\quad
 z[2] = \psi_6*\[A[1]^2*\psi_6 + B[1]*r'[1]]_{y[2]},\quad
 z[3] = \psi_6*\phi_6*C_5.
\label{eq.Bertini}
\]

Apart from the basepoint $(0:0:1)$ of the pencil, the fixed point locus of
this involution is the curve
\[
\deMaple
 K(y) := \psi_6*\[A[1]*y[3] + B[1]*y[1]]_{y[2]}
  -  \phi_6*\[A[2]*y[3] + B[3]*y[2]]_{y[1]} = 0.
\label{eq.K}
\]

\remark
The expressions for~$r$, $C_5$, and~$K$ found in~\cite{Moody} contain a
number of misprints.
The corrections suggested are verified by the identities in
\autoref{s.identities} below, as well as by~\eqref{eq.psi} and~\eqref{eq.K2}.
\endremark

\subsection{Further observations}\label{s.identities}
The expression for the Bertini involution, see~\eqref{eq.Bertini},
is obtained by substituting $z=lr+my$ and solving
$W_3(lr+my)=0$, see~\eqref{eq.W},
in $l:m$. (This equation is linear since $W_3(r)=W_3(y)=0$.)
Note that $\{\psi_6=0\}$ and $\{\phi_6=0\}$ are the curves contracted to the
basepoints $(1:0:0)$ and $(0:1:0)$, respectively. Hence, they can
also be found from the identities
\[*
\deMaple
y[3]*r[1] - y[1]*r[3] = A[2]*\gamma_4*\phi_6,\qquad
y[2]*r[3] - y[3]*r[2] = A[1]*\gamma_4*\psi_6.
\]
Besides, one has
\[*
\deMaple
y[1]*r[2] - y[2]*r[1] = -r'[1]*\gamma_4.
\]

A point~$y$ is fixed by~\eqref{eq.Bertini} if and only if the tangent at~$y$
to the member~\eqref{eq.pencil} of the pencil passing through~$y$ meets the
curve again at~$r$. In~\cite{Moody}, the equation~\eqref{eq.K} of the fixed
point locus is obtained by eliminating $\Gl:\mu$ from~\eqref{eq.pencil} and the
polar conic to~\eqref{eq.pencil} with respect to~$r$ (after the substitution
$x\mapsto y$).
Alternatively, $K$ can be found as the common factor of
$\deMaple y[3]*z[1]-y[1]*z[3]$ and $\deMaple y[2]*z[3]-y[3]*z[2]$,
using the identities
\[*
\deMaple
y[3]*z[1]-y[1]*z[3]=-\phi_6*K*A[2],\qquad
y[2]*z[3]-y[3]*z[2]=-\psi_6*K*A[1].
\]
Note that the rightmost factors are just two particular members of the
pencil.

\section{The map $\Cp2\da\Sigma_2$}\label{S.-2K}

From now on, we assume that the distinguished basepoint $(0:0:1)$ is simple.

\subsection{The anti-bicanonical map}\label{s.-2K}
Let~$Y$ be the plane~$\Cp2$ blown up at all basepoints
(including infinitely near) of the pencil other
than $(0:0:1)$. It is a (nodal, in general) Del Pezzo surface of degree~$1$,
and the anti-bicanonical linear system
maps $Y$ to a quadric cone in~$\Cp3$.
According to~\cite{Moody}, the proper transforms of the
sextics $\{\phi_6=0\}$ and $\{\psi_6=0\}$ are
in $\ls|-2K_Y|$. Hence, the space of sections $H^0(Y;-2K_Y)$ is generated by
$\phi_6$ (or $\psi_6$) and $w^2$, $ww'$, $w\2$, and the map
$y\mapsto\bz\in\Cp3$ is given by
\[*
\bz_0=\phi_6(y),\quad
\bz_1=w^2(y),\quad
\bz_2=w(y)w'(y),\quad
\bz_3=w\2(y).
\]
%where $(\bz_0:\bz_1:\bz_2:\bz_3)$ are homogeneous coordinates in~$\Cp3$.
Its image is the cone $\bz_1\bz_3=\bz_2^2$.
The passage to the affine coordinates $\bx:=\bz_1/\bz_2$, $\by:=\bz_0/\bz_2$
blows up the vertex and maps the cone to the Hirzebruch surface~$\Sigma_2$
with the exceptional section~$E$ of self-intersection~$(-2)$ (the exceptional
divisor over the vertex). The composed rational map $\Cp2\da\Sigma_2$ is
\[
\bx=w(y)/w'(y),\qquad
\by=\phi_6(y)/w\2(y).
\label{eq.map}
\]

Alternatively, $Y$ with the remaining basepoint $(0:0:1)$ blown up is a
rational Jacobian elliptic surface: the elliptic pencil is~\eqref{eq.pencil}
and the distinguished section is the exceptional divisor over~$(0:0:1)$.
The Bertini involution becomes the fiberwise multiplication by~$(-1)$, and
the quotient blows down to the Hirzebruch surface~$\Sigma_2$.
The quotient map is generically two-to-one; its ramification locus is the
union of the exceptional section $E\subset\Sigma_2$ and a certain proper
trigonal curve, \viz. the image of $\{K=0\}$.
The pull-backs of the fibers of~$\Sigma_2$ are the anti-canonical curves
in~$Y$ (\iq. the members of the original pencil~\eqref{eq.pencil} of cubics),
and the pull-backs of the \emph{proper} (\ie, disjoint from~$E$) sections
of~$\Sigma_2$ are the anti-bicanonical curves other than those representable
in the form $\{\Ga_1w^2+\Ga_2ww'+\Ga_3w\2=0\}$.

\subsection{The ramification locus}\label{s.ramification}
Since $\{\psi_6=0\}$ is the pull-back of a section of~$\Sigma_2$,
there must be a relation (after the substitution $x\mapsto y$) of the form
\[
\psi_6=\phi_6+S_2(w,w'),
\label{eq.psi}
\]
where $S_2$ is a certain homogeneous polynomial of degree~$2$, see below.

The curve $\{\phi_6=0\}$
is contracted by the Bertini
involution. Hence, the pull-back in~$Y$ of its image $\{\by=0\}$
splits into two components (sections of the elliptic pencil), of which one
is contracted by the blow down map $Y\to\Cp2$.
It follows that the free term~$R_3^2$ in equation~\eqref{eq.K2} below
is indeed a perfect square.

Since $\{K=0\}$ is the pull-back of the ramification locus (other
than~$E$), which is a proper trigonal curve,
there must be a relation
\[
K^2=-4\phi_6^3+\phi_6^2P_2(w,w')+\phi_6Q_4(w,w')+R_3^2(w,w'),
\label{eq.K2}
\]
where $P_2$, $Q_4$, and~$R_3$ are certain homogeneous polynomials of
degree~$2$, $4$, and~$3$, respectively. Let
$S_2(t_1,t_2)=\sum_{i=0}^2s_it_1^it_2^{2-i}$ \etc.
The coefficients $p_i$, $q_i$, $r_i$, $s_i$ are found by a direct
computation. They reduce to a remarkably simple form:
\[*
\deMaple
%\gathered
\aligned
 s[0] &= a[2]*c[1] - a[1]*c[2],\\
 r[0] &= -a[1]*b[2]*c[2]+a[1]*b[3]*c[1]+a[2]*b[1]*c[2],\\
 q[0] &= 4*(a[1]*c[2] - b[1]*b[3])*s[0]+2*b[2]*r[0],\\
 p[0] &= b[2]^2-4*a[2]*c[1]-4*b[1]*b[3]+8*a[1]*c[2]
\endaligned
\]
and
\[*
p_i=(-1)^i\{p_0\}_i,\quad
q_i=(-1)^i\{q_0\}_i,\quad
r_i=(-1)^i\{r_0\}_i,\quad
s_i=(-1)^i\{s_0\}_i,
%\endgathered
\]
where $\{\,\cdot\,\}_m$ is defined as follows: if $e$ is a
degree~$n$ monomial in $a_1,\ldots,c_2$, then $\{e\}_m$ is the sum of
$\binom{n}{m}$ monomials, each obtained from~$e$ by replacing $m$ of its~$n$
factors with their primed versions. (For example, one has
$\{a_1c_2\}_1=a_1c_2'+a_1'c_2$, $\{b_2^2\}_1=2b_2b_2'$,
and $\{a_1b_1c_1\}_2=a_1b_1'c_1'+a_1'b_1c_1'+a_1'b_1'c_1$.) This definition
extends to homogeneous polynomials by linearity.

\warning
The operation $\{\,\cdot\,\}_m$ is used only to shorten the notation. As with
the derivative, this operation should be performed \emph{before} any
substitution of any particular values of the coefficients (see, \eg, the
substitution $a_1=a_2=0$ in \autoref{S.Geiser}).
\endwarning

\remark\label{rem.symmetry}
Observe that
$S_2(w,w')$ remains unchanged under the transformation
$a_i\leftrightarrow a_i'$,
$b_i\leftrightarrow b_i'$,
$c_i\leftrightarrow c_i'$,
$w_i\leftrightarrow w_i'$.
The same holds for $P_2(w,w')$ and $Q_4(w,w')$, whereas $R_3(w,w')$ changes
sign.
\endremark

\problem
The symmetry in \autoref{rem.symmetry} is easily explained by
interchanging~$w$ and~$w'$. However, is there a geometric explanation for the
`regular' behaviour of the other coefficients?
\endproblem

Summarizing, we see that the map $\Cp2\da\Sigma_2$ given by~\eqref{eq.map}
takes the sextics $\{\phi_6=0\}$ and $\{\psi_6=0\}$ to the sections
$\{\by=0\}$ and $\{\by=-S_2(\bx)\}$, respectively.
The map is generically two-to-one (the deck translation being the Bertini
involution), and its ramification locus in~$\Sigma_2$ is the union of~$E$ and
the trigonal curve
\[*
-4\by^3+\by^2P_2(\bx)+\by Q_4(\bx)+R_3^2(\bx)=0.
\]
As usual, we treat the homogeneous bivariate polynomials $S_2$, $P_2$, $Q_4$,
and~$R_3$ as univariate ones
\via\ $S_2(\bx):=S_2(\bx,1)=\sum_{i=0}^2s_i\bx^i$ \etc.

\problem
Can one express coefficients~$s_i$ in terms of $p_i$, $q_i$, and~$r_i$? In
other words, does a choice of the ramification locus in~$\Sigma_2$ and one of
the sections select automatically the other section?
\endproblem

\subsection{Other sextics contracted by the involution}
The basepoint $(0:0:1)$ plays a special r\^ole in the definition of the
Bertini involution. The other basepoints are not special. In particular, any
other basepoint $(u_1:u_2:u_3)$ gives rise to a sextic $\{\psi^u_6=0\}$
contracted to this point and to a splitting section of~$\Sigma_2$ whose
pull-back this sextic is. Assuming that $u_1\ne0$ and normalizing the
coordinates as $(1,u_2,u_3)$,
we have
\[
\psi^u_6=\phi_6+S^u_2(w,w'),
\label{eq.other}
\]
where $S_2^u$ is a homogeneous polynomial of degree~$2$
whose coefficients $s^u_i$ are
\[*
\deMaple
s^u[0]=s[0]+(a[2]*c[2]*u_2+(a[2]*b[2]-a[1]*b[3])*u_3)+a[2]*b[3]*u_2*u_3+a[2]^2*u_3^2,
\qquad
s^u_i=(-1)^i\{s^u_0\}_i.
\]
As above, the image of $\{\psi^u_6=0\}$ in~$\Sigma_2$ is the section
$\{\by=-S^u_2(\bx)\}$.

These equations are easily obtained by changing the coordinates and placing
the basepoint in question to $(1:0:0)$.

\section{The Geiser involution}\label{S.Geiser}

Now, we assume that the pencil has exactly one basepoint infinitely near to
the distinguished point
$(0:0:1)$. In other words, all members of~\eqref{eq.pencil} have a common
tangent at $(0:0:1)$ and, hence, exactly one of them has a double point at
$(0:0:1)$. We can assume that this singular member is $\{w(x)=0\}$, thus
letting $a_1=a_2=0$.
The resulting special case of the Bertini involution is the
\emph{Geiser involution}, see~\cite{Moody}.

\subsection{The involution}\label{s.Geiser}
Most formulas in \autoref{s.Moody} simplify dramatically. One
obviously has
$\deMaple A[1] =-a'[1]*w(y)$ and $\deMaple A[2] =-a'[2]*w(y)$;
%\[*
%\deMaple
%A[1] =-a'[1]*w(y),\qquad
%A[2] =-a'[2]*w(y);
%\]
hence, $\gamma_4=w\gamma_1$, see~\eqref{eq.gamma}, where
\[
\gamma_1:=-(a'_1y_1+a'_2y_2)
\]
is the defining polynomial of the common tangent to the members of the pencil
at the distinguished basepoint $(0:0:1)$.
(Here and below, $w$ without an argument stands for $w(y)$.)

Next, there are splittings, see~\eqref{eq.R},
\[*
r'_1=w^2\tr'_1,\qquad
r'_3=w\tr'_3,\qquad
%r=(-a'_2\tr_1':a'_1\tr_1':a'_1a'_2\tr_3'),
r_i=w^3\tr_i,\ i=1,2,3,
\]
with
\[*
\gathered
\tr_1=-a'_2\tr_1',\qquad
\tr_2=a'_1\tr_1',\qquad
\tr_3=a'_1a'_2\tr_3',\\
\tr'_1=a\2_2B_1-a'_1a'_2B_2+a\2_1B_3,\qquad
\tr'_3=a'_1C_2-a'_2C_1.
\endgathered
\]
Furthermore, one has
\[*
\phi_6=w\phi_3,\qquad
\psi_6=w\psi_3,\qquad
C_5=w\tC,
\]
where
\[*
\deMaple
\aligned
\tC(y)&:=
 - a'[2]*\[B[1]-a'[1]*b[1]*y[1]*y[3]^2]_{y[2]}
 + a'[1]*\[a'[2]*y[3]*(w - b[1]*y[1]^2*y[3]) -  B[3]*y[2]]_{y[1]*y[2]},\\
\phi_3(y)&:=-a'[1]*C[2] + y[3]*\tC,\\
\psi_3(y)&:=-a'[2]*C[1] + y[3]*\tC.
\endaligned
\]
Finally, after reducing the common factor~$w^3$ in~\eqref{eq.Bertini},
the Geiser involution takes the form
\[*
\deMaple
 z[1] = \phi_3*\[a\2[2]*w*\phi_3 + B[3]*\tr'[1]]_{y[1]},\quad
 z[2] = \psi_3*\[a\2[1]*w*\psi_3 + B[1]*\tr'[1]]_{y[2]},\quad
 z[3] = \psi_3*\phi_3*\tC.
\label{eq.Bertini}
\]
The loci contracted to the
basepoints $(1:0:0)$ and $(0:1:0)$ are the cubics
$\{\psi_6=0\}$ and $\{\phi_6=0\}$, respectively, and the fixed point locus is
the sextic
\[*
\deMaple
 \tK(y) := \psi_3*\[-a'[1]*w*y[3] + B[1]*y[1]]_{y[2]}
  -  \phi_3*\[-a'[2]*w*y[3] + B[3]*y[2]]_{y[1]} = 0.
\]
One has $K=w\tK$, see~\eqref{eq.K}.
%Denoting $\gamma_1:=-(a'_1y_1+a'_2y_2)$ (\cf. the splitting
%of~$\gamma_4$ above),
The identities of \autoref{s.identities} turn into
\[*
\deMaple
\gathered
y[3]*\tr[1] - y[1]*\tr[3] = a'[2]*\gamma_1*\phi_3,\qquad
y[2]*\tr[3] - y[3]*\tr[2] = a'[1]*\gamma_1*\psi_3,\qquad
y[1]*\tr[2] - y[2]*\tr[1] = -\tr'[1]*\gamma_1,\\
y[3]*z[1]-y[1]*z[3]=a'[2]*\phi_3*\tK,\qquad
y[2]*z[3]-y[3]*z[2]=a'[1]*\psi_3*\tK.
\endgathered
\]

\subsection{The double covering $\Cp2\da\Cp2$}\label{s.quartic}
Let~$Z$ be the plane~$\Cp2$ blown up at the basepoints
(including infinitely near) of the pencil other
than $(0:0:1)$ and the infinitely near one.
It is a (nodal, in general) Del Pezzo surface of degree~$2$,
and the anti-canonical linear system
maps~$Z$ to~$\Cp2$. This map is generically two-to-one,
its deck translation is the Geiser involution, and its ramification
locus is a quartic curve in~$\Cp2$.
The anti-canonical system is the web of cubics passing through the seven
points blown up; the space $H^0(Z;-K_Z)$ is generated by any of~$\phi_3$
or~$\psi_3$ and by $w$ and~$w'$.
Hence, the corresponding rational map $\Cp2\da\Cp2$, $y\mapsto\bz$, is given
by
\[*
\bz_0=\phi_3(y),\quad
\bz_1=w(y),\quad
\bz_2=w'(y).
\]
It is straightforward that $q_0=r_0=s_0=0$. Hence, there are splittings
\[*
Q_4(t,t')=t\tQ_3(t,t'),\quad
R_3(t,t')=t\tR_2(t,t'),\quad
S_2(t,t')=t\tS_1(t,t')
\]
and relations~\eqref{eq.psi} and~\eqref{eq.K2} turn into
\[*
\gathered
\psi_3=\phi_3+\tS_1(w,w'),\\
\tK^2=-4\phi_3^3w+\phi_3^2P_2(w,w')+\phi_3\tQ_3(w,w')+\tR_2^2(w,w').
\endgathered
\]
Thus, the ramification locus is the quartic
\[
4\bz_0^3\bz_1=\bz_0^2P_2(\bz_1,\bz_2)+\bz_0\tQ_3(\bz_1,\bz_2)+\tR_2^2(\bz_1,\bz_2).
\label{eq.quartic}
\]
In the coordinates chosen, the lines $\{\bz_0=0\}$ and
$\{\bz_0+\tS_1(\bz_1,\bz_2)=0\}$ are double tangents (in the generalized sense)
to this quartic.

Similarly, given another basepoint $(1:u_2:u_3)$, one has
$S_2^u(t,t')=t\tS_1^u(t,t')$ and the cubic $\{\psi^u_3=0\}$ singular at this
point is given by (\cf.~\eqref{eq.other})
\[*
\psi_3^u=\phi_3+\tS_1^u(w,w').
\]

\remark
The coefficients of the polynomials $\tS_1$, $\tQ_3$, and~$\tR_2$ are the
same as those of~$S_2$, $Q_4$, and~$R_3$, respectively, see
\autoref{s.ramification}, with an obvious shift by one.
It is worth emphasizing that the brace operation $\{\,\cdot\,\}_i$ should be
evaluated \emph{before} the substitution $a_1=a_2=0$.
\endremark

\subsection{A few further observations}
Unlike the general case considered in \autoref{S.-2K},
now, the fixed point locus $\{\tK=0\}$ does pass through
%the distinguished basepoint
$(0:0:1)$. In fact, since this curve is the branch set, it follows
that $\{\tK=0\}$ is also the locus of the singular points of the singular
members of the web of cubics defined by the seven non-distinguished
basepoints. This curve has a double point at each of the seven basepoints. In
particular, it is an anti-bicanonical curve in~$Z$.

Under the natural identification of the web and the dual plane
$(\Cp2)\check{\,}$,
the locus of the singular members themselves is the curve dual
to~\eqref{eq.quartic}, including the lines through the singular points
of~\eqref{eq.quartic}.

As another observation, note that $\{\psi_3=0\}$ and $\{\phi_3=0\}$ are
special members of the web, \viz. those singular at $(1:0:0)$ and $(0:1:0)$,
respectively. As above, these cubics are contracted by the involution to the
corresponding basepoints.

\let\.\DOTaccent
\def\cprime{$'$}
\bibliographystyle{amsplain}
\bibliography{degt}

\providecommand{\bysame}{\leavevmode\hbox to3em{\hrulefill}\thinspace}
\providecommand{\MR}{\relax\ifhmode\unskip\space\fi MR }
\providecommand{\MRhref}[2]{%
  \href{http://www.ams.org/mathscinet-getitem?mr=#1}{#2}
}
\providecommand{\href}[2]{#2}
\begin{thebibliography}{1}

\bibitem{Moody}
Ethel~I. Moody, \emph{Notes on the {B}ertini involution}, Bull. Amer. Math.
  Soc. \textbf{49} (1943), 433--436. \MR{0008163 (4,253c)}

\end{thebibliography}

\end{document}